\documentclass{article}
\usepackage{amssymb}
\usepackage{amsfonts}
\usepackage{amsmath}

\newtheorem{theorem}{Theorem}

\newtheorem{corollary}[theorem]{Corollary}

\newtheorem{lemma}[theorem]{Lemma}

\newtheorem{proposition}[theorem]{Proposition}

\newenvironment{proof}[1][Proof]{\textbf{#1.} }{\ \rule{0.5em}{0.5em}}

\begin{document}

\title{\textsc{Hoeffding Decompositions and Two-Colour Urn Sequences}}
\author{Omar EL-DAKKAK\thanks{%
E-mail: omar.eldakkak@gmail.com} \ \ and \ \ Giovanni PECCATI\thanks{%
E-mail: giovanni.peccati@gmail.com} \\
\textsl{Laboratoire de Statistique Th\'{e}orique et Appliqu\'{e}e --
Universit\'{e} Paris VI \thanks{%
175, rue du Chevaleret, 75013 Paris, France.}}}
\date{October 19, 2006}
\maketitle

\begin{abstract}
Let $\mathbf{X}=\left( X_{1},X_{2},...\right) $ be a non-deterministic
infinite exchangeable sequence with values in $\left\{ 0,1\right\} $. We
show that $\mathbf{X}$ is Hoeffding-decomposable if, and only if, $\mathbf{X}
$ is either an i.i.d. sequence or a P\'{o}lya sequence. This completes the
results established in Peccati $\left[ 2004\right] $. The proof uses several
combinatorial implications of the correspondence between Hoeffding
decomposability and weak independence. Our results must be compared with
previous characterizations of i.i.d. and P\'{o}lya sequences given by Hill
\textit{et al. }$\left[ 1987\right] $ and Diaconis and Yilvisaker $\left[
1979\right] .$

\textbf{Key words. }Exchangeable Sequences; Hoeffding Decompositions; P\'{o}%
lya Urns; Weak Independence.

\textbf{AMS\ 2000 classifications: }60G09; 60G99
\end{abstract}

\section{Introduction, notation and statements\label{S : Intro}}

Let $\mathbf{X}_{\left[ 1,\infty \right) }$ $=$ $\{X_{n}$ $:$ $n\geq 1\}$ be
an exchangeable sequence of random observations, with values in some finite
set $D$. We say that $\mathbf{X}_{\left[ 1,\infty \right) }$ is \textsl{%
Hoeffding decomposable} if, for every $n\geq 2$, every symmetric statistic $%
T\left( X_{1},...,X_{n}\right) $ admits a unique representation as an
orthogonal sum of uncorrelated $U$-statistics with degenerate kernels of
increasing order. Hoeffding decompositions (also known as \textsl{ANOVA\
decompositions}) have been extensively studied for i.i.d. sequences (see for
instance Hoeffding [1948], Hajek [1968], Karlin and Rinott [1982], Vitale
[1991] and Bentkus \textit{et al.} [1997]) and for extractions without
replacement from a finite population (see Zhao and Chen [1990], Bloznelis
and G\"{o}tze [2001, 2002] and Bloznelis [2005]). In Peccati [2003, 2004,
2005], the second author of this paper has extended the theory of Hoeffding
decompositions to the framework of general exchangeable random sequences. In
particular, it was shown therein that the class of Hoeffding decomposable
exchangeable sequences coincides with the collection of weakly independent
sequences.

The aim of this note is to complete the results established in Peccati
[2004] by proving that \textsl{a (non deterministic) infinite exchangeable
sequence with values in }$\left\{ 0,1\right\} $\textsl{\ is Hoeffding
decomposable if, and only if, it is either a P\'{o}lya sequence or i.i.d.}.
As shown in Section \ref{S : FRem}, this result links the seemingly
unrelated notions of Hoeffding decomposable sequence and \textsl{urn process}%
, a concept studied e.g. in Hill \textsl{et al.} [1987].

Before stating our main theorem, we recall some basic definitions and facts
concerning Hoeffding decompositions and exchangeable sequences. We focus on
sequences with values in a finite set. The reader is referred to Peccati
[2004] for any unexplained concept or notation, as well as for general
statements concerning sequences with values in arbitrary Polish spaces.

\subsection{Preliminaries\label{SS : Pre}}

Let $D$ be a finite set, and consider an infinite exchangeable sequence $%
\mathbf{X}_{\left[ 1,\infty \right) }=\left\{ X_{n}:n\geq 1\right\} $ of $D$%
-valued random variables, defined on some probability space $\left( \Omega ,%
\mathcal{F},\mathbb{P}\right) $ such that $\mathcal{F}=\sigma (\mathbf{X}_{%
\left[ 1,\infty \right) })$. We recall that, according to the well-known
\textsl{de Finetti Theorem} (see e.g. Aldous [1983]), the assumption of
exchangeability implies that $\mathbf{X}_{\left[ 1,\infty \right) }$ is
necessarily a \textsl{mixture }of i.i.d. sequences with values in $D$.

For every $n\geq 1$ and every $1\leq u\leq n$, we write $\left[ n\right]
=\left\{ 1,...,n\right\} $ and $\left[ u,n\right] =\left\{
u,u+1,...,n\right\} $, and set $\mathbf{X}_{\left[ u,n\right] }\triangleq
\left( X_{u},X_{u+1}...,X_{n}\right) $ and $\mathbf{X}_{\left[ n\right]
}\triangleq \mathbf{X}_{\left[ 1,n\right] }=\left(
X_{1},X_{2},...,X_{n}\right) $. For every $n\geq 2$, we define the sequence
of spaces
\begin{equation*}
\left\{ SU_{k}\left( \mathbf{X}_{\left[ n\right] }\right) :k=0,...,n\right\}
,
\end{equation*}%
generated by symmetric $U$-statistics of increasing order, as follows: $%
SU_{0}\left( \mathbf{X}_{\left[ n\right] }\right) $ $\triangleq $ $\Re $
and, for $k=1,...,n$, $SU_{k}\left( \mathbf{X}_{\left[ n\right] }\right) $
is the collection of all random variables of the type%
\begin{equation}
F\left( \mathbf{X}_{\left[ n\right] }\right) =\sum_{1\leq
j_{1}<...<j_{k}\leq n}\varphi \left( X_{j_{1}},...,X_{j_{k}}\right) ,
\label{kUst}
\end{equation}%
where $\varphi $ is a real-valued symmetric function from $D^{k}$ to $\Re $.
A random variable such as $F$ in (\ref{kUst}) is called a $U$-\textsl{%
statistic with symmetric kernel of order} $k$. It is easily seen that the
kernel $\varphi $ appearing in (\ref{kUst}) is unique, in the sense that if $%
\varphi^{\prime }$ is another symmetric kernel satisfying (\ref{kUst}), then
$\varphi\left( \mathbf{X}_{\left[ k\right] }\right)=\varphi^{\prime }\left(
\mathbf{X}_{\left[ k\right] }\right)$, a.s.-$\mathbb{P}$. The following
facts are immediately checked: (i) for every $k=0,...,n$, $SU_{k}\left(
\mathbf{X}_{\left[ n\right] }\right) $ is a vector space, (ii) $%
SU_{k-1}\left( \mathbf{X}_{\left[ n\right] }\right) \subset SU_{k}\left(
\mathbf{X}_{\left[ n\right] }\right) $, (iii) $SU_{n}\left( \mathbf{X}_{%
\left[ n\right] }\right) =L_{s}\left( \mathbf{X}_{\left[ n\right] }\right) $%
, where (for $n\geq 1$) $L_{s}\left( \mathbf{X}_{\left[ n\right] }\right) $
is defined as the set of all random variables of the type $T\left( \mathbf{X}%
_{\left[ n\right] }\right) =T\left( X_{1},...,X_{n}\right) $, where $T$ is a
\textsl{symmetric function} from $D^{n}$ to $\Re $. The class of all
symmetric functions, from $D^{n}$ to $\Re $, will be denoted by $\mathcal{S}%
\left( D^{n}\right) $. Note that $L_{s}\left( \mathbf{X}_{\left[ n\right]
}\right) $ is a Hilbert space with respect to the inner product $%
\left\langle T_{1},T_{2}\right\rangle \triangleq \mathbb{E}\left[
T_{1}\left( \mathbf{X}_{\left[ n\right] }\right) T_{2}\left( \mathbf{X}_{%
\left[ n\right] }\right) \right] $, so that each $SU_{k}\left( \mathbf{X}_{%
\left[ n\right] }\right) $ is a closed subspace of $L_{s}\left( \mathbf{X}_{%
\left[ n\right] }\right) $. Finally, the sequence of \textsl{symmetric
Hoeffding spaces} $\left\{ SH_{k}\left( \mathbf{X}_{\left[ n\right] }\right)
:k=0,...,n\right\} $ associated to\ $\mathbf{X}_{\left[ n\right] }$ is
defined as $SH_{0}\left( \mathbf{X}_{\left[ n\right] }\right) $ $\triangleq $
$SU_{0}\left( \mathbf{X}_{\left[ n\right] }\right) $ $=$ $\Re $, and
\begin{equation}
SH_{k}\left( \mathbf{X}_{\left[ n\right] }\right) \triangleq SU_{k}\left(
\mathbf{X}_{\left[ n\right] }\right) \cap SU_{k-1}\left( \mathbf{X}_{\left[ n%
\right] }\right) ^{\perp }\text{, \ \ }k=1,...,n\text{,}  \label{k-hoeff}
\end{equation}%
where all orthogonals (here and in the sequel) are taken in $L_{s}\left(
\mathbf{X}_{\left[ n\right] }\right) $. Observe that $SH_{k}\left( \mathbf{X}%
_{\left[ n\right] }\right) \subset SU_{k}\left( \mathbf{X}_{\left[ n\right]
}\right) $ for every $k$, so that each $F\in SH_{k}\left( \mathbf{X}_{\left[
n\right] }\right) $ has necessarily the form (\ref{kUst}) for some
well-chosen symmetric kernel $\varphi $. Moreover, since $SU_{n}\left(
\mathbf{X}_{\left[ n\right] }\right) =L_{s}\left( \mathbf{X}_{\left[ n\right]
}\right) $, one has the following orthogonal decomposition:%
\begin{equation}
L_{s}\left( \mathbf{X}_{\left[ n\right] }\right)
=\bigoplus\limits_{k=0}^{n}SH_{k}\left( \mathbf{X}_{\left[ n\right] }\right)
\text{,}  \label{orto}
\end{equation}%
where \textquotedblleft\ $\oplus $ \textquotedblright\ stands for an
orthogonal sum. In particular, (\ref{orto}) implies that every symmetric
random variable $T\left( \mathbf{X}_{\left[ n\right] }\right) \in
L_{s}\left( \mathbf{X}_{\left[ n\right] }\right) $ admits a unique
representation as a non-correlated sum of $n+1$ terms, with the $k$th
summand ($k=0,...,n$) equal to an element of $SH_{k}\left( \mathbf{X}_{\left[
n\right] }\right) $.

The next definition, which is essentially borrowed from Peccati [2004],
formalizes the notion of \textquotedblleft Hoeffding
decomposability\textquotedblright\ evoked at the beginning of the section.

\bigskip

\textbf{Definition A.}\textsc{\ }The random sequence $\mathbf{X}_{\left[
1,\infty \right) }$ is \textsl{Hoeffding decomposable} if, for every $n\geq
2 $ and every $k=1,...,n$, the following double implication holds: $F\in
SH_{k}\left( \mathbf{X}_{\left[ n\right] }\right) $ if, and only if, the
kernel $\varphi $ appearing in its representation (\ref{kUst}) satisfies the
degeneracy condition
\begin{equation}
\mathbb{E}\left[ \varphi \left( \mathbf{X}_{\left[ k\right] }\right) \mid
\mathbf{X}_{\left[ 2,k\right] }\right] =0\text{, \ \ a.s.-}\mathbb{P}\text{.}
\label{degeneracy}
\end{equation}%
When a $U$-statistic $F$ as in (\ref{kUst}) is such that $\varphi $ verifies
(\ref{degeneracy}), one says that $F$ is a \textsl{completely degenerate}
symmetric $U$-statistic of order $k$, and that $\varphi $ is a \textsl{%
completely degenerate symmetric kernel} of order $k$.

For instance, when $k=3$, one has $\mathbf{X}_{\left[ 2,k\right] }=\left(
X_{2},X_{3}\right) $, and condition (\ref{degeneracy}) becomes:
\begin{equation*}
\mathbb{E}\left[ \varphi \left( X_{1},X_{2},X_{3}\right) \mid X_{2},X_{3}%
\right] =0.
\end{equation*}%
Of course, by exchangeability, (\ref{degeneracy}) holds if, and only if, $%
\mathbb{E}\left[ \varphi \left( \mathbf{X}_{\left[ k\right] }\right) \mid
\mathbf{X}_{\left[ k-1\right] }\right] $ $=$ $0$,\ a.s.-$\mathbb{P}$.

\bigskip

For every infinite non-deterministic exchangeable sequence $\mathbf{X}_{%
\left[ 1,\infty \right) }$ (not necessarily Hoeffding decomposable) and
every $k\geq 1$, the class of all kernels $\varphi :D^{k}\mapsto \Re $, such
that (\ref{degeneracy}) is verified, is noted $\Xi _{k}\left( \mathbf{X}_{%
\left[ 1,\infty \right) }\right)$.

\bigskip

It is well known (see e.g. Hoeffding [1948], Hajek [1968] or Karlin and
Rinott [1983]) that each i.i.d. sequence is decomposable in the sense of
Definition A. In Peccati [2004], the second author established a complete
characterization of Hoeffding decomposable sequences (with values in
arbitrary Polish spaces), in terms of \textsl{weak independence}. To
introduce this concept, we need some more notation. Fix $n\geq 2$, and
consider a symmetric function $T\in \mathcal{S}\left( D^{n}\right) $. We
define the function $\left[ T\right] _{n,n-1}^{\left( n-1\right) }$ as the
unique application from $D^{n-1}\ $to $\mathbb{\Re }$ such that%
\begin{equation}
\left[ T\right] _{n,n-1}^{\left( n-1\right) }\left( \mathbf{X}_{\left[ 2,n%
\right] }\right) =\mathbb{E}\left( T\left( \mathbf{X}_{\left[ n\right]
}\right) \mid \mathbf{X}_{\left[ 2,n\right] }\right) \text{, \ \ a.s.-}%
\mathbb{P}.  \label{Tn-1}
\end{equation}%
For instance, if $n=2$, then $\mathbf{X}_{\left[ 2\right] }=\left(
X_{1},X_{2}\right) $, $\mathbf{X}_{\left[ 2,2\right] }=X_{2}$ and $\left[ T%
\right] _{2,1}^{\left( 1\right) }\left( X_{2}\right) $ $=$ $\mathbb{E}\left(
T\left( X_{1},X_{2}\right) \mid X_{2}\right) $. Note that the
exchangeability assumption and the symmetry of $T$ imply that the
application $D^{n-1}\mapsto \mathbb{\Re }:\mathbf{x}\mapsto \left[ T\right]
_{n,n-1}^{\left( n-1\right) }\left( \mathbf{x}\right) $ is symmetric. Also,
with this notation, $T\in \Xi _{n}\left( \mathbf{X}_{\left[ 1,\infty \right)
}\right) $ if, and only if, $\left[ T\right] _{n,n-1}^{\left( n-1\right)
}\left( \mathbf{X}_{\left[ 2,n\right] }\right) =0$, a.s.-$\mathbb{P}$.

Analogously, for $u=2,...,n$ we define the function $\left[ T\right]
_{n,n-1}^{\left( n-u\right) }:D^{n-1}\mapsto \mathbb{\Re }$ through the
relation:
\begin{equation}
\left[ T\right] _{n,n-1}^{\left( n-u\right) }\left( \mathbf{X}_{\left[
u+1,u+n-1\right] }\right) =\mathbb{E}\left( T\left( \mathbf{X}_{\left[ n%
\right] }\right) \mid \mathbf{X}_{\left[ u+1,u+n-1\right] }\right) \text{, \
\ a.s.-}\mathbb{P}.  \label{Tn-u}
\end{equation}%
To understand our notation, observe that, for $u=2,...,n$, the two sets $%
\left[ n\right] $ and $\left[ u+1,u+n-1\right] $ have exactly $n-u$ elements
in common. For instance, if $n=3$ and $u=2$, then $\left[ u+1,u+n-1\right] $
$=$ $\left\{ 3,4\right\} $, and $\left[ T\right] _{3,2}^{\left( 1\right)
}\left( X_{3},X_{4}\right) $ $=$ $\mathbb{E}\left( T\left(
X_{1},X_{2},X_{3}\right) \mid X_{3},X_{4}\right) $. Again, exchangeability
and symmetry yield that the function $\mathbf{x}$ $\mapsto $ $\left[ T\right]
_{n,n-1}^{\left( 0\right) }\left( \mathbf{x}\right) $ (corresponding to the
case $u=n$) is symmetric on $D^{n-1}$. On the other hand, for $u=2,...,n-1$,
the application $\left( x_{1},...,x_{n-1}\right) $ $\mapsto $ $\left[ T%
\right] _{n,n-1}^{\left( n-u\right) }\left( x_{1},...,x_{n-1}\right) $ is
(separately) symmetric in the variables $\left( x_{1},...,x_{n-u}\right) $
and $\left( x_{n-u+1},...,x_{n-1}\right) $, and not necessarily symmetric as
a function on $D^{n-1}$. When $u=2,...,n-1$, we note $\widetilde{\left[ T%
\right] }_{n,n-1}^{\left( n-u\right) }$ the canonical symmetrization of $%
\left[ T\right] _{n,n-1}^{\left( n-u\right) }$.

Finally, for $u=2,...,n$, set
\begin{equation}
\widetilde{\Xi }_{n,n-u}\left( \mathbf{X}_{\left[ 1,\infty \right) }\right)
\triangleq \left\{ T\in \mathcal{S}\left( D^{n}\right) :\widetilde{\left[ T%
\right] }_{n,n-1}^{\left( n-u\right) }\left( \mathbf{X}_{\left[ u+1,u+n-1%
\right] }\right) =0\text{, \ \ a.s.-}\mathbb{P}\right\}  \label{phinu}
\end{equation}%
(recall that $\mathcal{S}\left( D^{n}\right) $ denotes the class of
symmetric functions on $D^{n}$). Note that, by exchangeability, $\widetilde{%
\left[ T\right] }_{n,n-1}^{\left( n-u\right) }\left( \mathbf{X}_{\left[
u+1,u+n-1\right] }\right) =0$,\ a.s.-$\mathbb{P}$, if, and only if, $%
\widetilde{\left[ T\right] }_{n,n-1}^{\left( n-u\right) }\left( \mathbf{X}_{%
\left[ n-1\right] }\right) =0$,\ a.s.-$\mathbb{P}$. The following technical
definition is taken from Peccati [2004].

\bigskip

\textbf{Definition B.}\textsc{\ }\label{Def : B}The exchangeable sequence $%
\mathbf{X}_{\left[ 1,\infty \right) }$ is \textsl{weakly independent }if,
for every $n\geq 2$,
\begin{equation}
\Xi _{n}\left( \mathbf{X}_{\left[ 1,\infty \right) }\right) \subset
\bigcap\limits_{u=2}^{n}\widetilde{\Xi }_{n,n-u}\left( \mathbf{X}_{\left[
1,\infty \right) }\right) \text{.}  \label{wi}
\end{equation}%
In other words, $\mathbf{X}_{\left[ 1,\infty \right) }$ is weakly
independent if, for every $n\geq 2$ and every $T\in \mathcal{S}\left(
D^{n}\right) $, the following implication holds: if $\left[ T\right]
_{n,n-1}^{\left( n-1\right) }\left( \mathbf{X}_{\left[ n-1\right] }\right)
=0 $, then $\widetilde{\left[ T\right] }_{n,n-1}^{\left( n-u\right) }\left(
\mathbf{X}_{\left[ n-1\right] }\right) =0$ for every $u=2,...,n$.

\bigskip

The next theorem, which is one of the main results of Peccati [2004], shows
that the notions of weak independence and Hoeffding decomposability are
equivalent for infinite exchangeable sequences.

\bigskip

\textsc{Theorem 0}\textbf{\ }(Peccati [2004, Th. 6])\textbf{. }\textit{%
Suppose that the infinite exchangeable sequence }$\mathbf{X}_{\left[
1,\infty \right) }$\textit{\ is such that, for every }$n\geq 2$\textit{, }%
\begin{equation}
SH_{k}\left( \mathbf{X}_{\left[ n\right] }\right) \neq \left\{ 0\right\}
\text{, \ \ }\forall k=1,...,n\text{.}  \label{non-t}
\end{equation}

\textit{Then, }$\mathbf{X}_{\left[ 1,\infty \right) }$\textit{\ is Hoeffding
decomposable if, and only if, it is weakly independent.}

\bigskip

\textbf{Remark. }Condition (\ref{non-t}) excludes for instance the case: $%
X_{n}=X_{1}$, for each $n\geq 1$.

\bigskip

Note that Theorem 0 also holds for exchangeable sequences with values in
general Polish spaces. In Peccati [2004] Theorem 0 has been used to show the
following two facts:

\begin{description}
\item[(F1)] There are infinite exchangeable sequences which are Hoeffding
decomposable and not i.i.d., as for instance the \textsl{Generalized Urn
Sequences} analyzed in Section 5 of Peccati [2004].

\item[(F2)] There exist infinite exchangeable sequences that \textsl{are not}
Hoeffding decomposable. For instance, one can consider a $\left\{
0,1\right\} $-valued exchangeable sequence $\mathbf{X}_{\left[ 1,\infty
\right) }^{Y}$ such that, conditioned on the realization of a random
variable $Y$ uniformly distributed on $\left( 0,\varepsilon \right) $ ($%
0<\varepsilon <1$), $\mathbf{X}_{\left[ 1,\infty \right) }^{Y}$ is composed
of independent Bernoulli trials with random parameter $Y$. See Peccati
[2004, p. 1807-1808] for more details.
\end{description}

Although the combination of Theorem 0, (F1) and (F2) gives several insights
into the structure of Hoeffding decomposable sequences, the analysis
contained in Peccati [2004] left open a crucial question: \textsl{can one
characterize the laws of Hoeffding decomposable sequences, in terms of their
de Finetti representation as mixtures of i.i.d. sequences?} In the following
sections, we will provide a complete answer when $D=\left\{ 0,1\right\} $,
by proving that in this case the class of Hoeffding decomposable sequences
contains exclusively i.i.d. and P\'{o}lya sequences. The extension of our
results to spaces $D$ with more than two elements is an open problem.

\subsection{Main results\label{SS : Main Res}}

For the rest of the paper, we will focus on the case $D=\left\{ 0,1\right\} $%
. According to the de Finetti Theorem, in this case the exchangeability of $%
\mathbf{X}_{\left[ 1,\infty \right) }=\left\{ X_{n}:n\geq 1\right\} $ yields
the existence of a probability measure $\gamma $ on $\left[ 0,1\right] $
such that, for every $n\geq 1$ and every vector $\left(
j_{1},...,j_{n}\right) \in \left\{ 0,1\right\} ^{n}$,%
\begin{equation}
\mathbb{P}\left\{ X_{1}=j_{1},...,X_{n}=j_{n}\right\} =\int_{\left[ 0,1%
\right] }\theta ^{\Sigma _{k}j_{k}}\left( 1-\theta \right) ^{n-\Sigma
_{k}j_{k}}\gamma \left( \mathrm{d}\theta \right) .  \label{miu}
\end{equation}

The measure $\gamma $,\ appearing in (\ref{miu}), is called the \textsl{de
Finetti measure} associated with $\mathbf{X}_{\left[ 1,\infty \right) }$. In
what follows, we shall systematically suppose that $\mathbf{X}_{\left[
1,\infty \right) }$ is \textsl{non-deterministic}, that is, that the support
of the measure $\gamma $ \textsl{is not contained} in $\left\{ 0\right\}
\cup \left\{ 1\right\} $. In particular, it is easy to prove that, when $%
D=\left\{ 0,1\right\} $, condition (\ref{non-t}) holds if, and only if, $%
\mathbf{X}_{\left[ 1,\infty \right) }$ is non-deterministic.

\bigskip

\textbf{Definition C.}\textsc{\ }The exchangeable sequence $\mathbf{X}_{%
\left[ 1,\infty \right) }=\left\{ X_{n}:n\geq 1\right\} $ is called a
\textsl{P\'{o}lya sequence} if there exist two real numbers $\alpha ,\beta
>0 $ such that
\begin{equation}
\gamma \left( \mathrm{d}\theta \right) =\frac{1}{B\left( \alpha ,\beta
\right) }\theta ^{\alpha -1}\left( 1-\theta \right) ^{\beta -1}\mathrm{d}%
\theta \text{,}  \label{beta}
\end{equation}%
where $\gamma $ is the de Finetti measure associated to $\mathbf{X}_{\left[
1,\infty \right) }$ through formula (\ref{miu}), and
\begin{equation*}
B\left( \alpha ,\beta \right) =\int_{0}^{1}\theta ^{\alpha -1}\left(
1-\theta \right) ^{\beta -1}\mathrm{d}\theta
\end{equation*}%
is the usual Beta function. The numbers $\alpha $ and $\beta $ are the
\textsl{parameters }of the P\'{o}lya sequence $\mathbf{X}_{\left[ 1,\infty
\right) }$. A random variable $\xi $, with values in $\left[ 0,1\right] $
and with law $\gamma $ as in (\ref{beta})$,$ is called a \textsl{Beta random
variable }of parameters $\alpha $ and $\beta $.

\bigskip

Classic references for the theory of P\'{o}lya sequences are Blackwell
[1973] and Blackwell and MacQueen [1973] (see also Pitman [1996, 2006] for a
state of the art review). Thanks to Peccati [2004, Corollary 9], we already
know that P\'{o}lya and i.i.d. sequences are Hoeffding decomposable. The
next result, which is the main achievement of our paper, shows that those
are the only exchangeable and Hoeffding decomposable sequences with values
in $\left\{ 0,1\right\} $. The proof is deferred to Section \ref{S : proofs}.

\begin{theorem}
\label{Th : Main}Let $\mathbf{X}_{\left[ 1,\infty \right) }$ be a
non-deterministic infinite exchangeable sequence of $\left\{ 0,1\right\} $%
-valued random variables. Then, the following two assertions are equivalent:

\begin{enumerate}
\item $\mathbf{X}_{\left[ 1,\infty \right) }$ is Hoeffding decomposable;

\item $\mathbf{X}_{\left[ 1,\infty \right) }$ is either an i.i.d. sequence
or a P\'{o}lya sequence.
\end{enumerate}
\end{theorem}

In Section \ref{S : FRem} we will discuss some connections between Theorem %
\ref{Th : Main} and the concept of \textsl{urn process,} as defined in Hill
\textit{et al.} [1987].

\bigskip

\textbf{Remarks. }We state two projection formulae, concerning respectively
i.i.d. and P\'{o}lya sequences.

(I) Let $\mathbf{X}_{\left[ 1,\infty \right) }$ be an i.i.d. sequence with
values in $\left\{ 0,1\right\} $, and fix $n\geq 2$ and $T\in L_{s}\left(
\mathbf{X}_{\left[ n\right] }\right) $. Then, for $k=1,...,n,$ the
projection of $T$ on the $k$-th Hoeffding space $SH_{k}\left( \mathbf{X}_{%
\left[ n\right] }\right) ,$ denoted by $\pi \left[ T,SH_{k}\right] ,$ is%
\begin{equation}
\pi \left[ T,SH_{k}\right] =\sum_{a=1}^{k}\left( -1\right) ^{k-a}\sum_{1\leq
j_{1}<\cdot \cdot \cdot <j_{a}\leq n}\left[ T-\mathbb{E}\left( T\right) %
\right] _{n,a}^{\left( a\right) }\left( X_{j_{1}},...,X_{j_{a}}\right) .
\label{sumaar}
\end{equation}
Formula (\ref{sumaar}) is classic (see e.g. Hoeffding [1948], or Vitale
[1991]), and can be easily deduced by an application of the
inclusion-exclusion principle.

(II) Let $\mathbf{X}_{\left[ 1,\infty \right) }$ be a P\'{o}lya sequence of
parameters $\alpha ,\beta >0$, and fix $n\geq 2$ and $T\in L_{s}\left(
\mathbf{X}_{\left[ n\right] }\right) $. Then, for $k=1,...,n,$ the
projection of $T$ on the $k$-th Hoeffding space associated with $\mathbf{X}_{%
\left[ n\right] }$ is of the form%
\begin{equation*}
\pi \left[ T,SH_{k}\right] =\sum_{a=1}^{k}\theta _{n}^{\left( k,a\right)
}\sum_{1\leq j_{1}<\cdot \cdot \cdot <j_{a}\leq n}\left[ T-\mathbb{E}\left(
T\right) \right] _{n,a}^{\left( a\right) }\left(
X_{j_{1}},...,X_{j_{a}}\right) .
\end{equation*}%
The explicit formulae describing the real coefficients $\theta _{n}^{\left(
k,a\right) }$ are given recursively in Peccati [2004, formula (23)]. For
instance, when $n=3,$ then%
\begin{equation*}
\left\{
\begin{array}{l}
\theta _{3}^{\left( 1,1\right) }=\frac{\alpha +\beta +1}{\alpha +\beta +2},
\\
\theta _{3}^{\left( 2,1\right) }=-\frac{\left( \alpha +\beta +1\right)
\left( \alpha +\beta +4\right) }{\left( \alpha +\beta +3\right) \left(
\alpha +\beta +2\right) }-\frac{\alpha +\beta +1}{\alpha +\beta +2}, \\
\theta _{3}^{\left( 2,2\right) }=\frac{\alpha +\beta +4}{\alpha +\beta +2}.%
\end{array}%
\right.
\end{equation*}

The rest of the paper is organized as follows: in Section \ref{S : Lemma} we
collect several technical results, leading to a new characterization of
Hoeffding decomposability in terms of conditional probabilities (see
Proposition \ref{P : ALtHD} below); the proof of Theorem \ref{Th : Main} is
contained in Section \ref{S : proofs}; in Section \ref{S : FRem}, a brief
discussion is presented, relating Theorem \ref{Th : Main} with several
notions associated with $\left\{ 0,1\right\} $-valued exchangeable sequences.

\section{Ancillary lemmas\label{S : Lemma}}

From now on, $\mathbf{X}_{\left[ 1,\infty \right) }=\left\{ X_{n}:n\geq
1\right\} $ will be a non-deterministic exchangeable sequence with values in
$D=\left\{ 0,1\right\} $. For $n\geq 2$, we write $\mathcal{S}\left( \left\{
0,1\right\} ^{n}\right) $ to indicate the vector space of symmetric
functions on $\left\{ 0,1\right\} ^{n}$, while $\mathfrak{S}_{n}$ stands for
the group of permutations of the set $\left[ n\right] =\left\{
1,...,n\right\} $. Given a vector $\mathbf{x}_{n}=\left(
x_{1},...,x_{n}\right) \in \left\{ 0,1\right\} ^{n}$ and a permutation $\pi
\in \mathfrak{S}_{n}$, we note $\mathbf{x}_{\pi \left( n\right) }$ the
action of $\pi $ on $\mathbf{x}_{n}$, that is, $\mathbf{x}_{\pi \left(
n\right) }=\left( x_{\pi \left( 1\right) },...,x_{\pi \left( n\right)
}\right) $. By exchangeability, we have of course that
\begin{equation*}
\mathbb{P}\left( \mathbf{X}_{\left[ n\right] }=\mathbf{x}_{n}\right) =%
\mathbb{P}\left( \mathbf{X}_{\left[ n\right] }=\mathbf{x}_{\pi \left(
n\right) }\right) \text{, \ \ }\forall n\geq 2\text{, }\forall \pi \in
\mathfrak{S}_{n}\text{,}
\end{equation*}%
yielding that, for $n\geq 2$, the value of the probability $\mathbb{P}\left(
\mathbf{X}_{\left[ n\right] }=\mathbf{x}_{n}\right) $ depends exclusively on
$n$ and on the number of zeros contained in the vector $\mathbf{x}_{n}$. For
$n\geq 1$ and $j=0,...,n$, we shall denote by $\mathbb{P}_{n}\left(
0^{\left( j\right) }\right) $ the \textsl{common value} taken by the
quantity $\mathbb{P}\left( \mathbf{X}_{\left[ n\right] }=\mathbf{x}%
_{n}\right) $ for all $\mathbf{x}_{n}=\left( x_{1},...,x_{n}\right) \in
\left\{ 0,1\right\} ^{n}$ such that $\mathbf{x}_{n}$ contains exactly $j$
zeros. For instance, when $n=3$ and $j=1$, one has that $\mathbb{P}%
_{3}\left( 0^{\left( 1\right) }\right) $ $=$ $\mathbb{P}\left( \mathbf{X}_{%
\left[ 3\right] }=\left( 0,1,1\right) \right) $ $=$ $\mathbb{P}\left(
\mathbf{X}_{\left[ 3\right] }=\left( 1,0,1\right) \right) $ $=$ $\mathbb{P}%
\left( \mathbf{X}_{\left[ 3\right] }=\left( 1,1,0\right) \right) $. Note
that, since $\mathbf{X}_{\left[ 1,\infty \right) }$ is non-deterministic, $%
\mathbb{P}_{n}\left( 0^{\left( j\right) }\right) >0$ for every $n\geq 1$ and
every $j=0,...,n$. Analogously, for every $n\geq 2$, every $j=0,...,n$, and
every symmetric function $\varphi \in \mathcal{S}\left( \left\{ 0,1\right\}
^{n}\right) $, we will write $\varphi \left( 0^{\left( j\right) }\right) $
to indicate the common value taken by $\varphi \left( \mathbf{x}_{n}\right) $
for all $\mathbf{x}_{n}\in \left\{ 0,1\right\} ^{n}$ containing exactly $j$
zeros.

The following result gives a complete characterization of the spaces%
\begin{equation*}
\Xi _{n}\left( \mathbf{X}_{\left[ 1,\infty \right) }\right) ,\qquad n\geq 2,
\end{equation*}%
defined through relation (\ref{degeneracy}) (note that, to define the spaces
$\Xi _{n}$ we do not need $\mathbf{X}_{\left[ 1,\infty \right)}$ to be
Hoeffding decomposable).

\begin{lemma}
\label{L : 1dimxin}With the assumptions and notation of this section, the
set $\Xi _{n}\left( \mathbf{X}_{\left[ 1,\infty \right) }\right) $ is the $1$%
-dimensional vector space spanned by the symmetric kernel $\varphi
_{n}^{\left( 0\right) }:\left\{ 0,1\right\} ^{n}\mapsto \Re $ defined by%
\begin{equation}
\varphi _{n}^{\left( 0\right) }\left( 0^{\left( k\right) }\right) =\left(
-1\right) ^{k}\frac{\mathbb{P}_{n}\left( 0^{\left( 0\right) }\right) }{%
\mathbb{P}_{n}\left( 0^{\left( k\right) }\right) },\qquad k=0,...,n.
\label{phinzero}
\end{equation}
\end{lemma}

\begin{proof}
Consider $\varphi _{n}\in \Xi _{n}\left( \mathbf{X}_{\left[ 1,\infty \right)
}\right) .$ By the definition of $\Xi _{n}\left( \mathbf{X}_{\left[ 1,\infty
\right) }\right) $, for any fixed $j=0,...,n-1$ and any fixed $\mathbf{x}%
_{n-1}\in \left\{ 0,1\right\} ^{n-1}$ such that $\sum_{i=1}^{n-1}\left(
1-x_{i}\right) =j,$ we have
\begin{eqnarray*}
0 &=&\mathbb{E}\left[ \varphi _{n}\left( \mathbf{X}_{\left[ n\right]
}\right) \mid \mathbf{X}_{\left[ 2,n\right] }=\mathbf{x}_{n-1}\right] \\
&=&\varphi _{n}\left( 0^{\left( j+1\right) }\right) \frac{\mathbb{P}%
_{n}\left( 0^{\left( j+1\right) }\right) }{\mathbb{P}_{n-1}\left( 0^{\left(
j\right) }\right) }+\varphi _{n}\left( 0^{\left( j\right) }\right) \frac{%
\mathbb{P}_{n}\left( 0^{\left( j\right) }\right) }{\mathbb{P}_{n-1}\left(
0^{\left( j\right) }\right) }\text{,}
\end{eqnarray*}%
and therefore $\varphi _{n}\left( 0^{\left( j+1\right) }\right) =-\left(
\mathbb{P}_{n}\left( 0^{\left( j\right) }\right) /\mathbb{P}_{n}\left(
0^{\left( j+1\right) }\right) \right) \times \varphi _{n}\left( 0^{\left(
j\right) }\right) .$ Arguing recursively on $j,$ one has%
\begin{equation}
\varphi _{n}\left( 0^{\left( j+1\right) }\right) =\left( -1\right) ^{j+1}%
\frac{\mathbb{P}_{n}\left( 0^{\left( 0\right) }\right) }{\mathbb{P}%
_{n}\left( 0^{\left( j+1\right) }\right) }\varphi _{n}\left( 0^{\left(
0\right) }\right) \text{, \ \ }j=0,...,n-1,  \label{recur}
\end{equation}%
showing that any symmetric kernel $\varphi _{n}\in \Xi _{n}\left( \mathbf{X}%
_{\left[ 1,\infty \right) }\right) $ is completely determined by the
quantity $\varphi _{n}\left( 0^{\left( 0\right) }\right) $. Now define a
kernel $\varphi _{n}^{\left( 0\right) }\in \Xi _{n}\left( \mathbf{X}_{\left[
1,\infty \right) }\right) $ by using (\ref{recur}) and by setting $\varphi
_{n}^{\left( 0\right) }\left( 0^{\left( 0\right) }\right) =\mathbb{P}%
_{n}\left( 0^{\left( 0\right) }\right) /\mathbb{P}_{n}\left( 0^{\left(
0\right) }\right) =1.$ It is easily seen that $\varphi _{n}^{\left( 0\right)
}$ must coincide with the function defined in (\ref{phinzero}). To conclude,
consider another element $\varphi _{n}$ of $\Xi _{n}\left( \mathbf{X}_{\left[%
1,\infty \right)}\right) $. Since there exists a constant $K\in \Re $ such
that $\varphi _{n}\left( 0^{\left( 0\right) }\right) =K=K\varphi
_{n}^{\left( 0\right) }\left( 0^{\left( 0\right) }\right) $, and since $%
\varphi _{n}$ has to satisfy (\ref{recur}), we deduce that $\varphi
_{n}=K\varphi _{n}^{\left( 0\right) }$, thus completing the proof.
\end{proof}

Given a function $f:\left\{ 0,1\right\} ^{m}\rightarrow \Re $, we note $%
\widetilde{f}$ its canonical symmetrization, that is: for every $\mathbf{x}%
_{m}\in \left\{ 0,1\right\} ^{m}$%
\begin{equation*}
\widetilde{f}\left( \mathbf{x}_{m}\right) =\frac{1}{m!}\sum_{\pi \in
\mathfrak{S}_{m}}f\left( \mathbf{x}_{\pi \left( m\right) }\right) \text{.}
\end{equation*}%
The following result will prove very useful.

\begin{lemma}
\label{L : Sym}Fix $m\geq 2$ and $v\in \left\{ 1,...,m-1\right\} $ and let
the application
\begin{equation*}
f_{v,m-v}:\left\{ 0,1\right\} ^{m}\mapsto \Re :\left( x_{1},...,x_{m}\right)
\mapsto f\left( x_{1},...,x_{m}\right) \text{,}
\end{equation*}%
be separately symmetric in the variables $\left( x_{1},...,x_{v}\right) $
and $\left( x_{v+1},...,x_{m}\right) $ (and not necessarily symmetric as a
function on $\left\{ 0,1\right\} ^{m}$). \textit{Then, for any} $\mathbf{x}%
_{m}=\left( x_{1},...,x_{m}\right) \in \left\{ 0,1\right\} ^{m}$ such that $%
\sum_{j=1}^{m}\left( 1-x_{j}\right) =z$ for some $z=0,...,m,$%
\begin{equation}
\widetilde{f}_{v,m-v}\left( \mathbf{x}_{m}\right) =\frac{\sum_{k=0\vee
\left( z-\left( m-v\right) \right) }^{z\wedge v}\binom{v}{k}\binom{m-v}{z-k}%
f_{v,m-v}\left( 0^{\left( k\right) },0^{\left( z-k\right) }\right) }{%
\sum_{k=0\vee \left( z-\left( m-v\right) \right) }^{z\wedge v}\binom{v}{k}%
\binom{m-v}{z-k}}.  \label{symNew}
\end{equation}%
where $f_{v,m-v}\left( 0^{\left( k\right) },0^{\left( z-k\right) }\right) $
denotes the common value of $f_{v,m-v}\left( \mathbf{y}_{m}\right) $ when $%
\mathbf{y}_{m}=\left( y_{1},...,y_{m}\right) $ is such that the vector $%
\left( y_{1},...,y_{v}\right) $ contains exactly $k$ zeros, and the vector $%
\left( y_{v+1},...,y_{m}\right) $ contains exactly $\left( z-k\right) $
zeros.

As a consequence, $\widetilde{f}_{v,m-v}\left( \mathbf{x}_{m}\right) =0$ for
every $\mathbf{x}_{m}\in \left\{ 0,1\right\} ^{m}$ if, and only if, for all $%
z=0,...,m$,%
\begin{equation}
\sum_{k=0\vee \left( z-\left( m-v\right) \right) }^{z\wedge v}\binom{v}{k}%
\binom{m-v}{z-k}f_{v,m-v}\left( 0^{\left( k\right) },0^{\left( z-k\right)
}\right) =0.  \label{Sym 0}
\end{equation}
\end{lemma}

\begin{proof}
Fix $\mathbf{x}_{m}\in \left\{ 0,1\right\} ^{m}$ such that $%
\sum_{j=1}^{m}\left( 1-x_{j}\right) =z$ for some $z=0,...,m.$ Without loss
of generality, we can assume%
\begin{equation*}
\mathbf{x}_{m}=(\underset{z\text{ times}}{\underbrace{0,0,...,0}},\underset{%
m-z\text{ times}}{\underbrace{1,1,...,1}}).
\end{equation*}%
Observe that, for all $k=\max \left\{ 0,z-\left( m-v\right) \right\}
,...,\min \left\{ z,v\right\} ,$ there are exactly $z!\left( m-z\right) !%
\binom{v}{k}\binom{m-v}{z-k}$ permutations $\pi \in \mathfrak{S}_{m}$ such
that $\sum_{j=1}^{v}\left( 1-x_{\pi \left( j\right) }\right) =k$ and $%
\sum_{j=v+1}^{m}\left( 1-x_{\pi \left( j\right) }\right) =z-k.$ The set of
all such permutations will be denoted by $\mathfrak{S}_{m}^{\left( k\right)
}.$ It is immediately seen that%
\begin{eqnarray*}
\widetilde{f}_{v,m-v}\left( \mathbf{x}_{m}\right) &=&\frac{1}{m!}%
\sum_{k=0\vee (z-\left( m-v\right))}^{z\wedge v}\sum_{\pi \in \mathfrak{S}%
_{m}^{\left( k\right) }}f_{v,m-v}\left( 0^{\left( k\right) },0^{\left(
z-k\right) }\right) \\
&=&\frac{1}{m!}\sum_{k=0\vee (z-\left( m-v\right))}^{z\wedge
v}f_{v,m-v}\left( 0^{\left( k\right) },0^{\left( z-k\right) }\right) \times
\mathrm{card}\left( \mathfrak{S}_{m}^{\left( k\right) }\right).
\end{eqnarray*}
Formula (\ref{symNew}) now follows by observing that%
\begin{equation*}
\frac{m!}{z!\left( m-z\right) !}=\binom{m}{z}=\sum_{k=0\vee (z-\left(
m-v\right)) }^{z\wedge v}\binom{v}{k}\binom{m-v}{z-k}.
\end{equation*}%
The last assertion in the statement of this lemma is an easy consequence of (%
\ref{symNew}).
\end{proof}

We shall conclude the section by obtaining a full characterization of $%
\left\{ 0,1\right\} $-valued Hoeffding decomposable sequences (stated in
Proposition \ref{P : ALtHD} below).

To do this, recall that, for any symmetric $\varphi :\left\{ 0,1\right\}
^{n}\mapsto \Re $, every $u=2,...,n$ and every $\mathbf{x}_{n-1}\in \left\{
0,1\right\} ^{n-1},$%
\begin{equation*}
\left[ \varphi \right] _{n,n-1}^{\left( n-u\right) }\left( \mathbf{x}%
_{n-1}\right) =\mathbb{E}\left( \varphi \left( \mathbf{X}_{\left[ n\right]
}\right) \mid \mathbf{X}_{\left[ u+1,u+n-1\right] }=\mathbf{x}_{n-1}\right) .
\end{equation*}%
Observe that the function $\left[ \varphi \right] _{n,n-1}^{\left(
n-u\right) }:$ $\left\{ 0,1\right\} ^{n-1}\mapsto \Re $ clearly meets the
symmetry properties of Lemma \ref{L : Sym} with $m=n-1$ and $v=n-u.$ Now fix
$z\in \left\{ 0,...,n-1\right\} $, and suppose that $\mathbf{x}_{n-1}\in
\left\{ 0,1\right\} ^{n-1}$ is such that $\sum_{j=1}^{n-1}\left(
1-x_{j}\right) =z$ and $\sum_{j=1}^{n-u}\left( 1-x_{j}\right) =k$. Then,%
\begin{equation}
\left[ \varphi \right] _{n,n-1}^{\left( n-u\right) }\left( \mathbf{x}%
_{n-1}\right) =\sum_{m=0}^{u}\binom{u}{m}\varphi \left( 0^{\left( k+m\right)
}\right) \frac{\mathbb{P}_{n-1+u}\left( 0^{\left( z+m\right) }\right) }{%
\mathbb{P}_{n-1}\left( 0^{\left( z\right) }\right) }.  \label{menta-1}
\end{equation}%
By applying (\ref{Sym 0}) in the case $m=n-1$ and $v=n-u$, we deduce that $%
\widetilde{\left[ \varphi \right] }_{n,n-1}^{\left( n-u\right) }\left(
0^{\left( z\right) }\right) =0$ if, and only if,%
\begin{equation}
\sum_{k=0\vee \left( z-\left( u-1\right) \right) }^{z\wedge \left(
n-u\right) }\binom{n-u}{k}\binom{u-1}{z-k}\left[ \varphi \right]
_{n,n-1}^{\left( n-u\right) }\left( 0^{\left( k\right) },0^{\left(
z-k\right) }\right) =0,  \label{menta}
\end{equation}%
where the notation $\widetilde{\left[ \varphi \right] }_{n,n-1}^{\left(
n-u\right) }\left( 0^{\left( z\right) }\right) $ and $\left[ \varphi \right]
_{n,n-1}^{\left( n-u\right) }\left( 0^{\left( k\right) },0^{\left(
z-k\right) }\right) $ has been introduced to indicate the value of $%
\widetilde{\left[ \varphi \right] }_{n,n-1}^{\left( n-u\right) }\left(
\mathbf{y}_{n-1}\right) $ (resp. $\left[ \varphi \right] _{n,n-1}^{\left(
n-u\right) }\left( \mathbf{w}_{n-1}\right) $), where $\mathbf{y}%
_{n-1}=\left( y_{1},...,y_{n-1}\right) \in \left\{ 0,1\right\} ^{n-1}$ is
any vector containing exactly $z$ zeros (resp. $\mathbf{w}_{n-1}=\left(
w_{1},...,w_{n-1}\right) \in \left\{ 0,1\right\} ^{n-1}$ is any vector
containing exactly $k$ zeros in $\left( w_{1},...,w_{n-u}\right) $ and $z-k$
zeros in $\left( w_{n-u+1},...,w_{n-1}\right) $).

Now recall that, by Theorem 0, $\mathbf{X}_{\left[ 1,\infty \right) }$ is
Hoeffding decomposable if, and only if, it is weakly independent, and that $%
\mathbf{X}_{\left[ 1,\infty \right) }$ is weakly independent if, and only
if, for all $n\geq 2$ and for any $\varphi \in \Xi _{n}\left( \mathbf{X}_{%
\left[ 1,\infty \right) }\right) $, one has $\varphi \in \widetilde{\Xi }%
_{n,u}\left( \mathbf{X}_{\left[ 1,\infty \right) }\right) $ for all $%
u=2,...,n.$ By Lemma \ref{L : 1dimxin}, we deduce that the sequence $\mathbf{%
X}_{\left[ 1,\infty \right) }$ is Hoeffding decomposable if, and only if,
for every $n\geq 2$ and every $u=2,...,n$, $\varphi _{n}^{\left( 0\right)
}\in \widetilde{\Xi }_{n,u}\left( \mathbf{X}_{\left[ 1,\infty \right)
}\right) $, where $\varphi _{n}^{\left( 0\right) }$ is defined in (\ref%
{phinzero}). By (\ref{menta}), this last relation is true if, and only if,
for every $n\geq 2$, every $z=0,...,n-1$ and every $u=2,...,n,$%
\begin{equation}
\sum_{k=0\vee \left( z-\left( u-1\right) \right) }^{z\wedge \left(
n-u\right) }\binom{n-u}{k}\binom{u-1}{z-k}\left[ \varphi _{n}^{\left(
0\right) }\right] _{n,n-1}^{\left( n-u\right) }\left( 0^{\left( k\right)
},0^{\left( z-k\right) }\right) =0.  \label{zero}
\end{equation}%
Substituting (\ref{phinzero}) and (\ref{menta-1}) in (\ref{menta}), we
obtain that (\ref{zero}) is true if, and only if,
\begin{eqnarray}
&&0=\frac{\mathbb{P}_{n}\left( 0^{\left( 0\right) }\right) }{\mathbb{P}%
_{n-1}\left( 0^{\left( z\right) }\right) }\sum_{k=0\vee \left( z-\left(
u-1\right) \right) }^{z\wedge \left( n-u\right) }\left( -1\right) ^{k}\binom{%
n-u}{k}\binom{u-1}{z-k}\times  \label{cinduriiin} \\
&&\text{ \ \ \ \ \ \ \ \ \ \ \ \ \ \ \ \ \ \ \ \ \ \ \ \ \ \ \ }\times
\sum_{m=0}^{u}\left( -1\right) ^{m}\binom{u}{m}\frac{\mathbb{P}%
_{n-1+u}\left( 0^{\left( m+z\right) }\right) }{\mathbb{P}_{n}\left(
0^{\left( m+k\right) }\right) }.  \notag
\end{eqnarray}

Note that%
\begin{equation}
\frac{\mathbb{P}_{n-1+u}\left( 0^{\left( m+z\right) }\right) }{\mathbb{P}%
_{n}\left( 0^{\left( m+k\right) }\right) }=\frac{1}{\binom{u-1}{z-k}}\mathbb{%
P}_{n+u-1}^{n}\left( 0^{\left( m+z\right) }\mid 0^{\left( m+k\right)
}\right) ,  \label{cindoro}
\end{equation}%
where $\mathbb{P}_{n+u-1}^{n}\left( 0^{\left( m+z\right) }\mid 0^{\left(
m+k\right) }\right) $ denotes the conditional probability that the vector $%
\mathbf{X}_{\left[ n+u-1\right] }$ contains exactly $m+z$ zeros, given that
the subvector $\mathbf{X}_{\left[ n\right] }$ contains exactly $m+k$ zeros.

\bigskip

\textbf{Remark. }For every $n\geq 1$, $0\leq a\leq b$, every $v\geq 1$, the
quantity $\mathbb{P}_{n+v}^{n}\left( 0^{\left( b\right) }\mid 0^{\left(
a\right) }\right) $ is equal to
\begin{equation*}
\mathbb{P}(\mathbf{X}_{\left[ n+1,n+v\right] }\text{ contains exactly }b-a%
\text{ zeros}\mid \mathbf{X}_{\left[ n\right] }\text{ contains exactly }a%
\text{ zeros}).
\end{equation*}

\bigskip

By plugging (\ref{cindoro}) into (\ref{cinduriiin}), we obtain the announced
characterization of weak independence.

\begin{proposition}
\label{P : ALtHD}Let $\mathbf{X}_{\left[ 1,\infty \right) }$ be a
non-deterministic infinite sequence of exchangeable $\left\{ 0,1\right\} $%
-valued random variables. For $\mathbf{X}_{\left[ 1,\infty \right) }$ to be
Hoeffding decomposable, it is necessary and sufficient that, for every $%
n\geq 2,$ every $u=2,...,n$ and every $z=0,...,n-1,$%
\begin{eqnarray}
&&0=\sum_{k=0\vee \left( z-\left( u-1\right) \right) }^{z\wedge \left(
n-u\right) }\left( -1\right) ^{k}\binom{n-u}{k}\times  \label{omaraus} \\
&&\text{ \ \ \ \ \ \ \ \ \ \ \ \ \ \ }\times \sum_{m=0}^{u}\left( -1\right)
^{m}\binom{u}{m}\mathbb{P}_{n+u-1}^{n}\left( 0^{\left( m+z\right) }\mid
0^{\left( m+k\right) }\right)  \notag
\end{eqnarray}
\end{proposition}

As shown in the next section, Proposition \ref{P : ALtHD} is the key tool to
prove Theorem \ref{Th : Main}.

\section{Proof of Theorem \protect\ref{Th : Main}\label{S : proofs}}

Here is an outline of the proof. We already know (thanks to Peccati [2004,
Corollary 9]) that, if $\mathbf{X}_{\left[ 1,\infty \right) }$ is either
i.i.d. or P\'{o}lya, then it is also Hoeffding decomposable, thus proving
the implication 2 $\Rightarrow $ 1. We shall therefore show that Hoeffding
decomposability implies necessarily that $\mathbf{X}_{\left[ 1,\infty
\right) }$ is either i.i.d. or P\'{o}lya$.$ The proof of this last
implication is divided in four steps. By using some easy remarks (Step 1)
and Proposition \ref{P : ALtHD}, we will prove that (\ref{omaraus}) implies
a universal relation linking the moments of the de Finetti measure $\gamma $
underlying any Hoeffding decomposable exchangeable sequence (Step 2). After
a discussion concerning the moments of Beta random variables (Step 3), we
conclude the proof in Step 4.

\bigskip

\textbf{Step 1. }We start with an easy remark. Define%
\begin{equation}
S\triangleq \left\{ \left( x,y,z\right) :0<x<y<z<1\right\} ,  \label{esset}
\end{equation}%
as well as the two functions%
\begin{eqnarray}
f\left( x,y,z\right) &=&2x^{2}z-xy^{2}-x^{2}y\text{, \ \ and}  \label{effe}
\\
g\left( x,y,z\right) &=&zx-2y^{2}+yz\text{.}  \label{ggi}
\end{eqnarray}%
\ Then, for any $\left( x,y,z\right) \in S$ one cannot have $f\left(
x,y,z\right) =0$ and $g\left( x,y,z\right) =0$ simultaneously.

\textbf{Step 2. }Let $\mathbf{X}_{\left[ 1,\infty \right) }=\left\{
X_{n}:n\geq 1\right\} $ be a non-deterministic exchangeable sequence with
values in $\left\{ 0,1\right\} $, and let $\gamma $ be the de Finetti
measure uniquely associated to $\mathbf{X}_{\left[ 1,\infty \right) }$
through formula (\ref{miu}). We denote by%
\begin{equation}
\mu _{n}=\mu _{n}\left( \gamma \right) =\int_{\left[ 0,1\right] }\theta
^{n}\gamma \left( \mathrm{d}\theta \right) \text{, \ \ }n\geq 0\text{,}
\label{mom}
\end{equation}%
the sequence of moments of $\gamma $ (the dependence on $\gamma $ is dropped
when there is no risk of confusion). We shall prove the following statement:
\textit{if }$\mathbf{X}_{\left[ 1,\infty \right) }$\textit{\ is Hoeffding
decomposable, then}%
\begin{equation}
\mu _{n+1}g\left( \mu _{n,}\mu _{n-1},\mu _{n-2}\right) =f\left( \mu
_{n},\mu _{n-1},\mu _{n-2}\right) \text{, \ \ }n\geq 2\text{,}  \label{mufg}
\end{equation}%
\textit{where }$f$\textit{\ and }$g$\textit{\ are respectively defined by (%
\ref{effe}) and (\ref{ggi}).}

To prove (\ref{mufg}), first recall that, due to Proposition \ref{P : ALtHD}%
, if $\mathbf{X}_{\left[ 1,\infty \right) }$ is Hoeffding decomposable, then
formula (\ref{omaraus}) must hold for every $n\geq 2$, every $u=2,...,n$ and
every $z=0,...,n-1$. In particular, it has to hold true for $u=2$, i.e., for
all $n\geq 2$ and all $z=0,...,n-1$, one must have that
\begin{equation}
\sum_{k=0\vee \left( z-1\right) }^{z\wedge \left( n-2\right) }\left(
-1\right) ^{k}\binom{n-2}{k}\sum_{m=0}^{2}\left( -1\right) ^{m}\binom{2}{m}%
\mathbb{P}_{n+1}^{n}\left( 0^{\left( m+z\right) }\mid 0^{\left( m+k\right)
}\right) =0,  \label{u2}
\end{equation}%
for every $n\geq 2$ and every $z=0,...,n-1$. For $z=0$, formula (\ref{u2})
becomes%
\begin{equation}
\mathbb{P}_{n+1}^{n}\left( 0^{\left( 2\right) }\mid 0^{\left( 2\right)
}\right) -2\mathbb{P}_{n+1}^{n}\left( 0^{\left( 1\right) }\mid 0^{\left(
1\right) }\right) +\mathbb{P}_{n+1}^{n}\left( 0^{\left( 0\right) }\mid
0^{\left( 0\right) }\right) =0.  \label{z0}
\end{equation}%
For $z=n-1,$ (\ref{u2}) is equivalent to%
\begin{equation}
\mathbb{P}_{n+1}^{n}\left( 0^{\left( n\right) }\mid 0^{\left( n\right)
}\right) -2\mathbb{P}_{n+1}^{n}\left( 0^{\left( n-1\right) }\mid 0^{\left(
n-1\right) }\right) +\mathbb{P}_{n+1}^{n}\left( 0^{\left( n-2\right) }\mid
0^{\left( n-2\right) }\right) =0.  \label{zn-1}
\end{equation}%
For $1\leq z\leq n-2,$ (\ref{u2}) becomes%
\begin{eqnarray}
0 &=&\binom{n-2}{z-1}[\mathbb{P}_{n+1}^{n}\left( 0^{\left( z+2\right) }\mid
0^{\left( z+1\right) }\right)  \label{ZZZ} \\
&&\text{ \ \ \ \ \ \ \ \ \ \ \ \ \ \ \ \ }-2\mathbb{P}_{n+1}^{n}\left(
0^{\left( z+1\right) }\mid 0^{\left( z\right) }\right) +\mathbb{P}%
_{n+1}^{n}\left( 0^{\left( z\right) }\mid 0^{\left( z-1\right) }\right) ]
\notag \\
&&-\binom{n-2}{z}[\mathbb{P}_{n+1}^{n}\left( 0^{\left( z+2\right) }\mid
0^{\left( z+2\right) }\right)  \notag \\
&&\text{ \ \ \ \ \ \ \ \ \ \ \ \ \ \ \ \ }-2\mathbb{P}_{n+1}^{n}\left(
0^{\left( z+1\right) }\mid 0^{\left( z+1\right) }\right) +\mathbb{P}%
_{n+1}^{n}\left( 0^{\left( z\right) }\mid 0^{\left( z\right) }\right) ].
\notag
\end{eqnarray}%
Combining (\ref{z0}), (\ref{zn-1}) and (\ref{ZZZ}), we deduce that (\ref{u2}%
) is true if, and only if, for all $p=0,...,n-2,$%
\begin{equation}
\mathbb{P}_{n+1}^{n}\left( 0^{\left( p+2\right) }\mid 0^{\left( p+2\right)
}\right) -2\mathbb{P}_{n+1}^{n}\left( 0^{\left( p+1\right) }\mid 0^{\left(
p+1\right) }\right) +\mathbb{P}_{n+1}^{n}\left( 0^{\left( p\right) }\mid
0^{\left( p\right) }\right) =0.  \label{zp}
\end{equation}%
Now, for $p=0,...,n-2$, write $\Delta _{p}$ to indicate the (forward)
difference operator of order $p$, given by: $\Delta _{0}f\left( n\right)
=f\left( n\right) $, $\Delta _{1}f\left( n\right) =f\left( n+1\right)
-f\left( n\right) ,$ and
\begin{equation*}
\Delta _{p}=\underset{p\text{ times}}{\underbrace{\Delta _{1}\circ \cdot
\cdot \cdot \circ \Delta _{1}}}\text{.}
\end{equation*}%
For $p=0,...,n-2$ one has%
\begin{eqnarray*}
\mathbb{P}_{n+1}^{n}\left( 0^{\left( p\right) }\mid 0^{\left( p\right)
}\right) &=&\mathbb{P}_{n+1}^{n}\left( X_{n+1}=1\mid \mathbf{X}_{\left[ n%
\right] }\text{ contains }p\text{ zeros}\right) \\
&=&\frac{\Delta _{p}\mu _{n+1-p}}{\Delta _{p}\mu _{n-p}}\text{,}
\end{eqnarray*}%
where the sequence of moments $\mu _{n}$, $n\geq 1$, is given by (\ref{mom}%
). Since (\ref{zp}) must hold for $p=0,$ we deduce that%
\begin{equation*}
\frac{\Delta _{2}\mu _{n-1}}{\Delta _{2}\mu _{n-2}}-2\frac{\Delta _{1}\mu
_{n}}{\Delta _{1}\mu _{n-1}}+\frac{\mu _{n+1}}{\mu _{n}}=0,
\end{equation*}%
and straightforward calculations yield relation (\ref{mufg}).

\bigskip

\textbf{Remark. }Suppose that $\mathbf{X}_{\left[ 1,\infty \right) }$ is
exchangeable and non-deterministic, and define $\mu _{n}$, $n\geq 0$, via (%
\ref{mom}). Then, we have that $\mu _{n+1}\in \left( 0,1\right) $ for every $%
n\geq 0$, and that, for every $n\geq 2$, $\left( \mu _{n},\mu _{n-1},\mu
_{n-2}\right) \in S$, where $S$ is defined as in (\ref{esset}). As a
consequence, the conclusions of Step 1 and (\ref{mufg}) imply that, if $%
\mathbf{X}_{\left[ 1,\infty \right) }$ is Hoeffding decomposable, then $%
f\left( \mu _{n},\mu _{n-1},\mu _{n-2}\right) \neq 0$ and $g\left( \mu
_{n},\mu _{n-1},\mu _{n-2}\right) \neq 0$ for every $n\geq 2$. Therefore,

\begin{equation}
\mu _{n+1}=\frac{f\left( \mu _{n},\mu _{n-1},\mu _{n-2}\right) }{g\left( \mu
_{n},\mu _{n-1},\mu _{n-2}\right) }.  \label{mufg2}
\end{equation}

\bigskip

\textbf{Step 3. }We claim that, for any $\left( c_{1},c_{2}\right) \in
\left( 0,1\right) ^{2}$ such that $c_{1}^{2}<c_{2}<c_{1}$, there exists a
unique pair $\left( \alpha ^{\ast },\beta ^{\ast }\right) \in \left(
0,+\infty \right) \times \left( 0,+\infty \right) $ such that
\begin{equation*}
\mathbb{E}\left[ \xi \right] =c_{1}\text{ \ \ and \ \ }\mathbb{E}\left[ \xi
^{2}\right] =c_{2},
\end{equation*}%
where $\xi $ is a Beta random variable of parameters $\alpha ^{\ast }$ and $%
\beta ^{\ast }$. To check this, just observe that, if $\xi $ is Beta of
parameters $\alpha $ and $\beta $ then%
\begin{equation*}
\mathbb{E}\left( \xi \right) =\frac{\alpha }{\alpha +\beta }\text{ \ \ and \
\ }\mathbb{E}\left( \xi ^{2}\right) =\frac{\alpha \left( \alpha +1\right) }{%
\left( \alpha +\beta \right) \left( \alpha +\beta +1\right) }\text{, }
\end{equation*}%
and that, for every fixed $\left( c_{1},c_{2}\right) \in \left( 0,1\right)
^{2}$ such that $c_{1}^{2}<c_{2}<c_{1}$, the system%
\begin{equation}
\left\{
\begin{array}{l}
\frac{\alpha }{\alpha +\beta }=c_{1} \\
\frac{\alpha \left( \alpha +1\right) }{\left( \alpha +\beta \right) \left(
\alpha +\beta +1\right) }=c_{2}%
\end{array}%
\right. ,  \label{foma}
\end{equation}%
$\allowbreak $admits a unique solution $\left( \alpha ^{\ast },\beta ^{\ast
}\right) \in \left( 0,+\infty \right) \times \left( 0,+\infty \right) $:
namely%
\begin{equation}
\left\{
\begin{array}{l}
\alpha ^{\ast }=\frac{c_{1}\left( 1-c_{2}\right) }{c_{2}-c_{1}^{2}}, \\
\beta ^{\ast }=\frac{\left( 1-c_{1}\right) \left( c_{1}-c_{2}\right) }{%
c_{2}-c_{1}^{2}}.%
\end{array}%
\right.  \label{systemoon}
\end{equation}

We are now in a position to conclude the proof of the implication 1 $%
\Rightarrow $ 2 in the statement of Theorem \ref{Th : Main}.

\textbf{Step 4. }Let $\mathbf{X}_{\left[ 1,\infty \right) }$ be a
non-deterministic exchangeable sequence, denote by $\gamma $ its de Finetti
measure and by $\left\{ \mu _{n}\left( \gamma \right) :n\geq 0\right\} $ the
sequence of moments appearing in (\ref{mom}). We suppose that $\mathbf{X}_{%
\left[ 1,\infty \right) }$ is Hoeffding decomposable. There are only two
possible cases: either $\mu _{1}\left( \gamma \right) ^{2}=\mu _{2}\left(
\gamma \right) $, or $\mu _{1}\left( \gamma \right) ^{2}<\mu _{2}\left(
\gamma \right) $. If $\mu _{1}\left( \gamma \right) ^{2}=\mu _{2}\left(
\gamma \right) $, then necessarily $\gamma =\delta _{x}$ for some $x\in
\left( 0,1\right) $, and therefore $\mathbf{X}_{\left[ 1,\infty \right) }$
is a sequence of i.i.d. Bernoulli trials with common parameter equal to $x$.
If $\mu _{1}\left( \gamma \right) ^{2}<\mu _{2}\left( \gamma \right) $,
then, thanks to the results contained in Step 3 (note that $\mu _{2}\left(
\gamma \right) <\mu _{1}\left( \gamma \right) $, since $\mathbf{X}_{\left[
1,\infty \right) }$ is non-deterministic), there exists a unique pair $%
\left( \alpha ^{\ast },\beta ^{\ast }\right) \in \left( 0,+\infty \right)
\times \left( 0,+\infty \right) $ such that
\begin{eqnarray}
\mu _{1}\left( \gamma \right) &=&\mathbb{E}\left( \xi \right) =\frac{1}{%
B\left( \alpha ^{\ast },\beta ^{\ast }\right) }\int_{0}^{1}\theta \theta
^{\alpha ^{\ast }-1}\left( 1-\theta \right) ^{\beta ^{\ast }-1}\mathrm{d}%
\theta  \label{a} \\
\mu _{2}\left( \gamma \right) &=&\mathbb{E}\left( \xi ^{2}\right) =\frac{1}{%
B\left( \alpha ^{\ast },\beta ^{\ast }\right) }\int_{0}^{1}\theta ^{2}\theta
^{\alpha ^{\ast }-1}\left( 1-\theta \right) ^{\beta ^{\ast }-1}\mathrm{d}%
\theta \text{,}  \label{aa}
\end{eqnarray}%
where $\xi $ stands for a Beta random variable of parameters $\alpha ^{\ast
} $ and $\beta ^{\ast }$. Moreover, (\ref{mufg2}) and the fact that P\'{o}%
lya sequences are Hoeffding decomposable imply that, for any $n\geq 2$,%
\begin{eqnarray*}
\mu _{n+1}\left( \gamma \right) &=&\frac{f\left( \mu _{n}\left( \gamma
\right) ,\mu _{n-1}\left( \gamma \right) ,\mu _{n-2}\left( \gamma \right)
\right) }{g\left( \mu _{n}\left( \gamma \right) ,\mu _{n-1}\left( \gamma
\right) ,\mu _{n-2}\left( \gamma \right) \right) }\text{, \ \ and } \\
\mathbb{E}\left( \xi ^{n+1}\right) &=&\frac{f\left( \mathbb{E}\left( \xi
^{n}\right) ,\mathbb{E}\left( \xi ^{n-1}\right) ,\mathbb{E}\left( \xi
^{n-2}\right) \right) }{g\left( \mathbb{E}\left( \xi ^{n}\right) ,\mathbb{E}%
\left( \xi ^{n-1}\right) ,\mathbb{E}\left( \xi ^{n-2}\right) \right) }\text{,%
}
\end{eqnarray*}%
where $f$ and $g$ are given by (\ref{effe}) and (\ref{ggi}). As (\ref{a})
and (\ref{aa}) are in order, we deduce that, for every $n\geq 1$,
\begin{equation}
\mu _{n}\left( \gamma \right) =\mathbb{E}\left( \xi ^{n}\right) =\frac{1}{%
B\left( \alpha ^{\ast },\beta ^{\ast }\right) }\int_{0}^{1}\theta ^{n}\theta
^{\alpha ^{\ast }-1}\left( 1-\theta \right) ^{\beta ^{\ast }-1}\mathrm{d}%
\theta \text{.}  \label{aaa}
\end{equation}%
Since probability measures on $\left[ 0,1\right] $ are determined by their
moments, the combination of (\ref{a}), (\ref{aa}) and (\ref{aaa}) gives
\begin{equation*}
\gamma \left( \mathrm{d}\theta \right) =\frac{1}{B\left( \alpha ^{\ast
},\beta ^{\ast }\right) }\theta ^{\alpha ^{\ast }-1}\left( 1-\theta \right)
^{\beta ^{\ast }-1}\mathrm{d}\theta \text{,}
\end{equation*}%
implying that $\mathbf{X}_{\left[ 1,\infty \right) }$ is a P\'{o}lya
sequence of parameters $\alpha ^{\ast }$ and $\beta ^{\ast }$. This
concludes the proof of Theorem \ref{Th : Main}.

\section{Further remarks\label{S : FRem}}

\textbf{(I) }With the terminology of Hill \textit{et al. }[1987], a random
sequence $\mathbf{X}_{\left[ 1,\infty \right) }=\left\{ X_{n}:n\geq
1\right\} $, with values in $\left\{ 0,1\right\} $, is called an \textsl{urn
process }if there exists a measurable function $f:\left[ 0,1\right] $ $%
\mapsto $ $\left[ 0,1\right] $ and positive natural numbers $r,b>0$, such
that, for every $n\geq 1$,
\begin{equation}
\mathbb{P}\left( X_{n+1}=1\mid X_{1},...,X_{n}\right) =f\left( \frac{%
r+X_{1}+\cdot \cdot \cdot +X_{n}}{r+b+n}\right) .  \label{urns}
\end{equation}%
According to Theorem 1 in Hill \textit{et al. }[1987], the only exchangeable
and non-deterministic urn processes are i.i.d. and P\'{o}lya sequences with
integer parameters (for which $f$ is, respectively, constant and equal to
the identity map). This yields immediately the following consequence of
Theorem \ref{Th : Main}, showing that the two (seemingly unrelated) notions
of urn process and Hoeffding decomposable sequence are in many cases
equivalent. The proof can be achieved by using the calculations performed in
Step 4.

\begin{corollary}
\label{C : urn}Let $\mathbf{X}_{\left[ 1,\infty \right) }=\left\{
X_{n}:n\geq 1\right\} $ be a $\left\{ 0,1\right\} $-valued infinite
exchangeable non-deterministic sequence such that%
\begin{equation}
\mathbb{P}\left( X_{1}=1\right) =c_{1}\qquad \text{and\qquad }\mathbb{P}%
\left( X_{1}=X_{2}=1\right) =c_{2},  \label{matsuri}
\end{equation}%
for some constants $c_{1}$ and $c_{2}$ such that $0<c_{1}^{2}<c_{2}<c_{1}<1.$
If the system (\ref{foma}) admits integer solutions, then $\mathbf{X}_{\left[
1,\infty \right) }$ is Hoeffding decomposable if, and only if, it is an urn
process.

In general, a sequence $\mathbf{X}_{\left[ 1,\infty \right) }$ verifying (%
\ref{matsuri}) is Hoeffding decomposable if, and only if, it is a P\'{o}lya
sequence with parameters $\alpha ^{\ast }$ and $\beta ^{\ast }$ given by (%
\ref{systemoon}).
\end{corollary}

\bigskip

\textbf{(II) }The arguments rehearsed in the proof of Theorem \ref{Th : Main}
provide an alternative proof of Theorem 5 in Diaconis and Ylvisaker [1979].
Indeed, in this reference it is shown that, if an exchangeable sequence is
such that its predictive probabilities
\begin{equation*}
\mathbb{P}_{n+1}^{n}\left( 0^{\left( p\right) }\mid 0^{\left( p\right)
}\right) =\mathbb{P}\left( X_{n+1}=1\mid \mathbf{X}_{\left[ n\right] }\text{
contains exactly }p\text{ zeros}\right)
\end{equation*}%
depend linearly on $p$, then its de Finetti measure must be\ Beta or Dirac.
To see how Diaconis and Ylvisaker's result can be recovered using our
techniques, suppose that a given exchangeable random sequence is not i.i.d.
and is such that its predictive probabilities verify the equation%
\begin{equation}
\mathbb{P}_{n+1}^{n}\left( 0^{\left( p\right) }\mid 0^{\left( p\right)
}\right) =a_{n}p+b_{n}  \label{questa}
\end{equation}%
for some positive sequences $\left\{ a_{n}\right\} $ and $\left\{
b_{n}\right\} $. Then, it is immediately seen that $\mathbb{P}%
_{n+1}^{n}\left( 0^{\left( p\right) }\mid 0^{\left( p\right) }\right) $ also
verifies (\ref{zp}), and one deduces from the previous discussion that the
associated de Finetti measure must be Beta. Conversely, if one supposes that
the predictive probabilities $\mathbb{P}_{n+1}^{n}\left( 0^{\left( p\right)
}\mid 0^{\left( p\right) }\right) $ of an exchangeable non-i.i.d. sequence
verify the difference equation (\ref{zp}), then one must conclude that the $%
\mathbb{P}_{n+1}^{n}\left( 0^{\left( p\right) }\mid 0^{\left( p\right)
}\right) $ verifies (\ref{questa}) and that
\begin{equation*}
a_{n}=\frac{1}{1+a\left( n-1\right) }\text{ \ and \ }b_{n}=\frac{b}{%
1+a\left( n-1\right) },
\end{equation*}%
for some $a>0,$ $b>0$ such that $a+b<1.$

\bigskip

\textbf{Acknowledgement. }Both authors wish to thank Igor Pr\"{u}nster for
inspiring discussions as this work progressed.

\end{document}